# COMPARING FINITE TYPE INVARIANTS OF KNOTS AND INTEGRAL HOMOLOGY 3-SPHERES

STAVROS GAROUFALIDIS

This edition: February 6, 1996     First edition: February 1,1996

Fax number: (401) 863 9013

Email: `stavros@math.mit.edu`

ABSTRACT. Using elementary counting methods of weight systems for finite type invariants of knots and integral homology 3-spheres, in the spirit of [B-NG], we answer positively three questions raised in [Ga]. In particular, we exhibit a one-to-one map from the space of finite type invariants of integral homology 3-spheres to the space of finite type invariants of knots in $S^3$.

## Contents

## 1. INTRODUCTION

**1.1. History.** Finite type invariants of integral homology 3-spheres have been very recently introduced by T. Ohtsuki in his seminal paper [Oh]. Subsequently, and similar to the theory of finite type invariants of knots (see [B-N] and [Ko]) weight

The author was partially supported by NSF grant DMS-95-05105. This and related preprints can also be obtained by accessing the WEB in the address `http://www.math.brown.edu/~stavrosg/`

systems for 3-manifold invariants were introduced, and recently the universal finite type invariant of integral homology 3-spheres was constructed [LMO], and [L]. The purspose of the present paper is to use elementary counting arguments of weight systems in the spirit of [B-N] to answer in positive three questions raised in [Ga].

**1.2. Preliminaries.** We begin by recalling some standard termonilogy from the theory of finite type invariants of knots, see [B-N]. let $\mathcal{K}$ denote the vector space over $\mathbb{Q}$ on the set of isotopy classes of oriented knots in $S^3$. For a nonegative integer $m$, let $\mathcal{F}_m\mathcal{K}$ denote the subspace of $\mathcal{K}$ spanned by all knots with $m$ double points, according to figure 1. It is easy to see that this defines a decreasing filtration on $\mathcal{K}$. A ($\mathbb{Q}$-linear) function $v : \mathcal{K} \to \mathbb{Q}$ is called a *type $m$* invariant of knots (in $S^3$) if $v(\mathcal{F}_{m+1}\mathcal{K}) = 0$. Let $\mathcal{F}_m\mathcal{V}$ denote the space of type $m$ invariants of knots. For $v$ a finite type invariant of knots, let $W(v)$ denote the associated *weight system*, see [B-N]. The notion of weight systems is a combinatorial notion that reduces questions on knot invariants to questions on (finite, but usually exponentially growing) combinatorics.

**Figure 1.** A figure resolving double points and defining the filtration on $\mathcal{K}$.

We also need to recall from [Oh] and [Ga] some standard terminology from the theory of finite type invariants of oriented integral homology 3-spheres. All 3-manifolds will be smooth and oriented and all links will be oriented. Following [Oh] and [Ga], let $\mathcal{M}$ denote the vector space over $\mathbb{Q}$ on the set of orientation preserving diffeomorphism classes of integral homology 3-spheres. A link $L$ in an integral homology 3-sphere is called *AS-admissible* if the linking numbers between its components vanish, and if every component is assigned an element of $\{-1, +1\}$ (such an assignement will be called a *unit framing*, and will generically be denoted by $f$). For an $AS$-admissible link $(L, f)$ in an integral homology 3-sphere $M$, let $M_{L,f}$ denote the manifold obtained by Dehn surgery on the framed link $L$. For such a link $(L, f)$, let

$$[M, L, f] = \sum_{L' \subseteq L} (-1)^{|L'|} M_{L',f'} \in \mathcal{M} \tag{1}$$

where the sum is over all sublinks of $L$ (including the empty one), $f'$ is the restriction of the framing $f$ of $L$ to $L'$ and $|L'|$ is the number of components of $L'$. Note that since $(L, f)$ is an $AS$-admissible link, $M_{(L',f')}$ is a integral homology 3-sphere for every sublink $L'$ of $L$. Let $\mathcal{F}_m\mathcal{M}$ denote the subspace of $\mathcal{M}$ spanned by all elements of the form $[M, L, f]$ for all $AS$-admissible $m$-component links in integral homology 3-spheres $M$. It is easy to see that this defines a decreasing filtration on $\mathcal{M}$. We call a map $\mathcal{M} \to \mathbb{Q}$ a *type $m$* invariant of integral homology 3-spheres if $v(\mathcal{F}_{m+1}\mathcal{M}) = 0$. Let $\mathcal{F}_m\mathcal{O}$ denote the space of type $m$ invariants of integral homology 3-spheres. In [GO] we introduced the notion of a *weight system* $W(v)$ for a finite type invariant $v$ of integral homology 3-spheres.

In [Ga] we considered the map: $\mathcal{K} \to \mathcal{M}$ defined by $K \to [S^3, K, +1] = S^3 - S^3_{K,+1}$. This is a classical map, used often in the study of knots (or 3-manifolds). The above map dually induces one:

$$Hom(\mathcal{K}, \mathbb{Q}) \to Hom(\mathcal{M}, \mathbb{Q}) \tag{2}$$

It is trivial to show that [Ga] the above map maps $\mathcal{F}_m\mathcal{O} \to \mathcal{F}_{m-1}\mathcal{V}$. Furthermore, in [Ga] we asked the following questions:

**Q1:** Does the above map factor through a map:

$$\Phi_m : \mathcal{F}_{3m}\mathcal{O} \to \mathcal{F}_{2m}\mathcal{V} \tag{3}$$

**Q2:** If it does, is the map $\Phi_m$ one-to-one?

**Q3:** With the terminology of [B-NG], is it true that for $v \in \mathcal{F}_m\mathcal{O}$, the associated weight system $W(\Phi_m(v))$ is a product of the coefficients of the Alexander-Conway weight system?

Since then, progress has been made to the above questions. N. Habegger [Ha] and independently the author in joint work with J. Levine [GL1] have positively answered question 1. The purspose of the present paper is to answer in positive all the above questions, using the methods of weight systems. This in particular exhibits the relation between the *methods* of [B-N], [KSA] and [LMO].

**1.3. Statement of the results.** We are now ready to state the results of the paper. Let $m$ be an arbitraty nonnegative integer.

**Theorem 1.** [Ha][GL1] *The map of equation (2) factors through a map $\Phi_m$ of equation (3).*

**Theorem 2.** *The map $\Phi_m$ is one-to-one.*

**Theorem 3.** *For $v \in \mathcal{F}_m\mathcal{O}$, the associated weight system $W(\Phi_m(v))$ is a product of the coefficients of the Alexander-Conway weight system.*

A few remarks are in order:

*Remark* 1.1. As mentioned before, the method of proving theorems 1, 2 and 3 is the use of weight systems, which reduces a topological problem to a combinatorial one. This method, used for the first time in [B-NG] blends well with the results of [KSA], [LMO] and [L].

*Remark* 1.2. Note that since the result of Dehn surgery is independent of the orientation of the knot, it implies that the image $\Phi_m(\mathcal{F}_{3m}\mathcal{O})$ are finite type invariants of oriented knots that cannot detect the change of the orientation of the knot.

*Remark* 1.3. Note that for $m \neq 1$, $v \in \mathcal{F}_{3m}\mathcal{O}$ then $\Phi_m(v)$ is *not* a canonical finite type invariant of knots, with the terminology of [B-NG]. It can be expanded in terms canonical invariants of knots (of degree at most $2m$). Studying this expansion, may shed some light in the subdiagonal terms of the expansion of the colored Jones polynomial (see [B-N]). It may also shed some light in subleading terms of finite type invariants of integral homology 3-spheres (see [GL2]).

*Remark* 1.4. In [GL2] it was shown that every integral homology 3-sphere can be obtained by surgery on a unit-framed boundary link in $S^3$. As a corollary, we obtained the curious corollary that the Casson invariant of integral homology 3-spheres can be calculated in terms of the Casson invariant for knots. Theorem 2 is a generalization of this corollary at all degrees (the case $m = 1$ being the Casson invariant).

*Remark* 1.5. Recall that there are two main notions that produce invariants of knots and 3-manifolds: the notion of *topological quantum field theory* ($TQFT$ for short) and the notion of finite type. If we compose the map (2) with a $TQFT$, we can ask an analogous set of questions. For example, question 2 asks whether for a $TQFT$, the (numeraical) knot invariants determine the 3-manifold invariants. The answer seems to be positive, with the appropriate nondegeneracy assumptions, the reason being that a $TQFT$ comes equipped with finite dimensional representations of the braid group, and knowing the invariants for knots determines the invariants for links, and thus for 3-manifolds.

*Remark* 1.6. The map $\Phi_m$ is by no means onto. For example, for $m = 2$ we have that $dim\mathcal{F}_6\mathcal{O} = 4$ but $dim\mathcal{F}_4\mathcal{K} = 6$.

**1.4. Acknowledgment.** We wish to thank T. Le and T. Ohtsuki for useful electronic conversations, and the `Internet` for providing a means for such conversations.

## 2. Proofs

In this section we recall a few fundamental facts about finite type invariants of knots and integral homology 3-spheres which will be used to prove theorems 1, 2 and 3.

**2.1. The universal finite type invariant of knots and integral homology 3-spheres.** The present section is devoted to a review of the universal finite type invariant of knots and integral homology 3-spheres. The material here has been drown from [Ko], [B-N], [LMO] and [L].

We begin by recalling that a *UT-graph* is a graph with univalent and trivalent vertices only, together with an orientation (i.e., cyclic order) to each trivalent vertex. For an arbitrary compact, connected, oriented, 1-dimensional manifold $X$ (possibly the empty one!) a *chord diagram* with support on $X$ is the manifold $X$ together with a $UT$-graph whose univalent vertices are on $X$ and such that the graph has no components isomorphic to the sircle. Two chord diagrams with support $X$ are *equal* if there is an orientation preserving diffeomorphism of the $UT$-graphs that extends to one on $X$. The *degree* of a chord diagram with support on $X$ is half the number of vertices of the associated $UT$-graph. Thus a degree $m$ chord diagram with support the empty set is a $UT$-graph without univalent vertices, with $3m$ edges and $2m$ vertices. Let $\mathcal{A}(X)$ denote the vector space (over $\mathbb{Q}$) on the set of equivalent classes of chord diagrams with support $X$ modulo two relations: the *antisymmetry AS* and the *IHX*, see figure 2. Note that our version of the $IHX$ relation includes the $STU$ relation of [LMO]. Note also that $\mathcal{A}(X)$ is a graded vector space. Let $\mathcal{A}_m(X)$ denote the degree $m$ subspace of $\mathcal{A}(X)$. In the present paper we will be interested in the particular cases of $X = S^1$ and $X = \emptyset$. In both of these cases, $\mathcal{A}(X)$ has a a multiplication (see [LMO]), which makes $\mathcal{A}(X)$ a commutative graded algebra. Let $\hat{\mathcal{A}}(X)$ denote the completed algebra.

**Figure 2.** The $AS$ and the $IHX$ relations on $\mathcal{A}(X)$.

With these preliminaries in mind, we can state the fundamental theorem of finite type invariants of knots and integral homology 3-spheres, due to [Ko] and [LMO], [L]:

**Theorem 4.** *For every nonnegative integer $m$ we have the following short exact sequences:*

$$1 \longrightarrow \mathcal{F}_{m+1}\mathcal{K} \longrightarrow \mathcal{F}_m\mathcal{K} \xrightarrow{W_m^K} \mathcal{A}_m(S^1) \longrightarrow 1$$

$$1 \longrightarrow \mathcal{F}_{3m+3}\mathcal{M} \longrightarrow \mathcal{F}_m\mathcal{M} \xrightarrow{W_m^M} \mathcal{A}_m(\emptyset) \longrightarrow 1$$

*where $W_m^K$ is the map defined in [Ko], [B-N] (and was shown there to be onto), and $W_m^M$ is the map defined in [GO] (and was shown in [LMO] and [L] to be onto).*

In particular, one obtains the following corollary:

**Corollary 2.1.** *The above maps induce isomorphisms of finite dimensional vector spaces as follows:*

(4) $\quad\quad \mathcal{F}_m\mathcal{K} = Hom(\mathcal{A}_{\leq m}(S^1), \mathbb{Q}) \quad and \quad \mathcal{F}_{3m}\mathcal{M} = Hom(\mathcal{A}_{\leq m}(\emptyset), \mathbb{Q})$

The maps considered above (depending on $m$) assemble well to define two maps:

(5) $\quad\quad\quad\quad Z^K : \mathcal{K} \to \hat{\mathcal{A}}(S^1) \ \text{ and } \ Z^{LMO} : \mathcal{M} \to \hat{\mathcal{A}}(\emptyset)$

The image under these maps of a knot or an integral homology 3-sphere satisfies a *universal property* (due to [Ko] in the case of knots, and due to [L] in the case of integral homology 3-spheres) and is therefore referred to as the *universal finite type invariant* of knots and integral homology 3-spheres respectively. These two universal finite type invariants of knots and integral homology 3-spheres are related to each other. In fact, recalling the map $\mathcal{K} \to \mathcal{M}$ from the introduction, we have the following commutative diagram:

$$\begin{array}{ccc} \mathcal{K} & \xrightarrow{[S^3,,+1]} & \mathcal{M} \\ Z^K \downarrow & & \downarrow Z^{LMO} \\ \hat{\mathcal{A}}(S^1) & \xrightarrow{\bar{\iota}} & \hat{\mathcal{A}}(\emptyset) \\ \downarrow & & \downarrow \\ \mathcal{A}_{\leq 2m}(S^1) & \xrightarrow{\bar{\iota}_m} & \mathcal{A}_{\leq m}(\emptyset) \end{array}$$

where the map $\bar{\iota}$ is constucted as the inverse limit of the maps $\bar{\iota}_m$, and the map $\bar{\iota}_m$ is defined (following the notation of [LMO]) in such a way that makes the following diagram commutative:

$$\begin{array}{ccc} \mathcal{A}_{\leq 2m}(S^1) & \xrightarrow{\bar{\iota}_m} & \mathcal{A}_{\leq m}(\emptyset) \\ \downarrow & & \downarrow \\ \mathcal{A}^{\circ}_{\leq 2m}(S^1)/L_{<2m}P_{m+1}O_m & \xrightarrow{\iota_m} & \mathcal{A}^{\circ}_{\leq m}(\emptyset)/L_{<2m}P_{m+1}O_m \end{array}$$

The notation has been borrowed by [LMO]. $\mathcal{A}^{\circ}(X)$ has been defined as in [LMO, section 3.1], the relations $L_{<2m}P_{m+1}O_m$ as in [LMO, proposition 3.3], and $\iota_m$ (the main combinatorial map) as in [LMO, section 2]. In [LMO, lemma 3.1] it was shown that the first vertical map in the above diagram is onto, and in [LMO, lemma 3.4] that the second once is an isomorphism therefore implying the existence of the map $\bar{\iota}_m$.

**2.2. The proofs.** Collecting the above information, we get the following commutative diagram:

$$\begin{array}{ccc} \mathcal{K} & \xrightarrow{[S^3,,+1]} & \mathcal{M} \\ \downarrow & & \downarrow \\ \mathcal{A}_{\leq 2m}(S^1) & \xrightarrow{\bar{\iota}_m} & \mathcal{A}_{\leq m}(\emptyset) \end{array}$$

where the vertical maps and $\bar{\iota}_m$ are onto. This together with corollary 2.1 proves that the map (2) factors through a map $\Phi_m$ of equation 3, thus proving theorem 1. Furthermore, since the vertical maps on the above diagram and $\bar{\iota}_m$ are onto are onto, this shows that dually the map $\Phi_m$ is one-to-one, thus showing theorem 2.

In order to show theorem 3, let $v \in \mathcal{F}_{3m}\mathcal{O}$ be a type $3m$ invariant of integral homology 3-spheres, let $\Phi_m(v) \in \mathcal{F}_{2m}\mathcal{V}$ be the associated type $2m$ knot invariant and let $W(\Phi_m(v)) \in Hom(\mathcal{A}_{2m}(S^1), \mathbb{Q})$ be the associated weight of the knot invariant.

With the definition of [B-N, section 6], we can recall the following result from [KSA]:

**Theorem 5.** [KSA] *A weight system $W \in Hom(\mathcal{A}_{2m}(S^1), \mathbb{Q})$ is the product of the coefficients of the Conway-Alexander polynomial if and only if it vanishes for all chord diagrams whose associated UT-graph has less than $2m$ univalent vertices.*

Now we can give a proof of theorem 3 as follows: $v \in \mathcal{F}_{3m}\mathcal{O}$ defines a map $W(v) \in Hom(\mathcal{A}_m(\emptyset), \mathbb{Q})$ and by composition with $\overline{\iota}_m$ a map in $Hom(\mathcal{A}_{2m}(S^1), \mathbb{Q})$. Due to the commutativity of the above diagram this map is the knot weight system $W(\Phi_m(v))$. On the other hand, for every chord diagram that satisfies the hypothesis of theorem 5 has zero image in $\mathcal{A}^\circ_{\leq 2m}(S^1)/L_{<2m}P_{m+1}O_m$ (due to relation $L_{<2m}$), and therefore, due to the above commutative diagrams, we have $W(\Phi_m(v)) = 0$. This together with theorem 5 finishes the proof of theorem 3.

151 THAYER STREET, DEPARTMENT OF MATHEMATICS, BROWN UNIVERSITY PO BOX 1917, PROVIDENCE RI 02912-0001. U.S.A.

*E-mail address*: stavros@math.mit.edu